\documentclass[12pt]{amsart}

\usepackage[dvips]{graphicx}
\usepackage{float}
\usepackage[below]{placeins}
\usepackage{amssymb, amsthm} 
\usepackage{hyperref}
\newtheorem{theorem}{Theorem}

\newtheorem{claim}[theorem]{Claim}

\newtheorem{lemma}[theorem]{Lemma}

\theoremstyle{remark}

\newcommand{\R}{\mathbf{R}}
\newcommand{\pd}{ {\partial}}

\newcommand{\eps}{\varepsilon}

\newcommand{\htau}{^{(h, \tau)}}

\begin{document}

\title[Construction of Self-Similar Surfaces under MCF]{Construction of Complete Embedded Self-Similar Surfaces under Mean Curvature Flow. Part I. }
\author{Xuan Hien Nguyen}
\address{Department of Mathematics, University of Wisconsin - Madison, WI 53706} 
\curraddr{Department of Mathematical Sciences, University of Cincinnati, PO Box 210025, Cincinnati, OH 45221}
\email{hien.nguyen@uc.edu}

\begin{abstract}
We carry out the first main step towards the construction of new examples of complete embedded  self-similar surfaces under mean curvature flow. An approximate solution is obtained by taking two known examples of self-similar surfaces and  desingularizing the intersection circle  using an appropriately modified singly periodic Scherk surface, called the core. Using an inverse function theorem, we show that for small boundary conditions on the core, there is an embedded surface close to the core that is a solution of the equation for self-similar surfaces. This provides us with an adequate central piece to substitute for  the intersection.
\end{abstract}

\maketitle

\section{Introduction}

This paper is the first one of a series of three articles describing the construction of new examples of complete embedded self-similar surfaces under mean curvature flow \cite{mine;part2}\cite{mine;part3}.  Our general strategy is inspired by  Kapouleas'  article  \cite{kapouleas;embedded-minimal-surfaces}. His success in constructing complete embedded minimal surfaces motivates us to adapt  his method towards finding self-similar surfaces under mean curvature flow.

Self-similar solutions  are solutions to the mean curvature flow that do not change shape but are merely contracted (called self-shrinkers) or dilated (self-expanders) by it. Up to scaling, the self-shrinkers satisfy the equation
	\begin{equation}
	\label{eq:self-shrinker}
	H+X\cdot \nu =0,
	\end{equation}
where $H$ is the mean curvature and $\nu$ is the normal vector so that the mean curvature vector is $\mathbf{H}=H \nu$. The sign of $H$ is chosen so that the mean curvature of a convex surface is positive. These solutions are of special interest because they model the behavior of the mean curvature flow at singularities if the blow-up is of type I.  There are currently only four known embedded complete self-shrinkers: a plane, a sphere, a cylinder and a shrinking doughnut \cite{angenent;doughnuts}.  The availability of new examples will lead to a better understanding of the behavior of the flow near singularities. 

\subsection{Strategy}

\begin{figure}[b]
\centering
\includegraphics[height=2.3in]{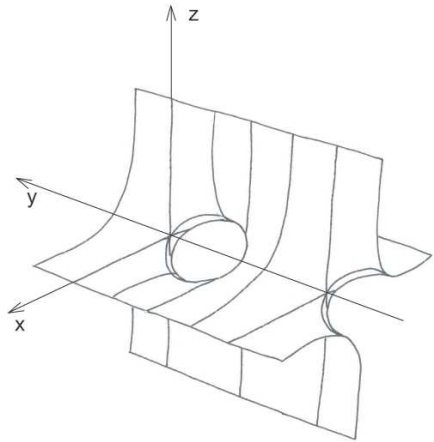}
\caption{Scherk surface}
\label{fig:scherk}
\end{figure}

\begin{figure}[t]

\includegraphics[height=1.6in]{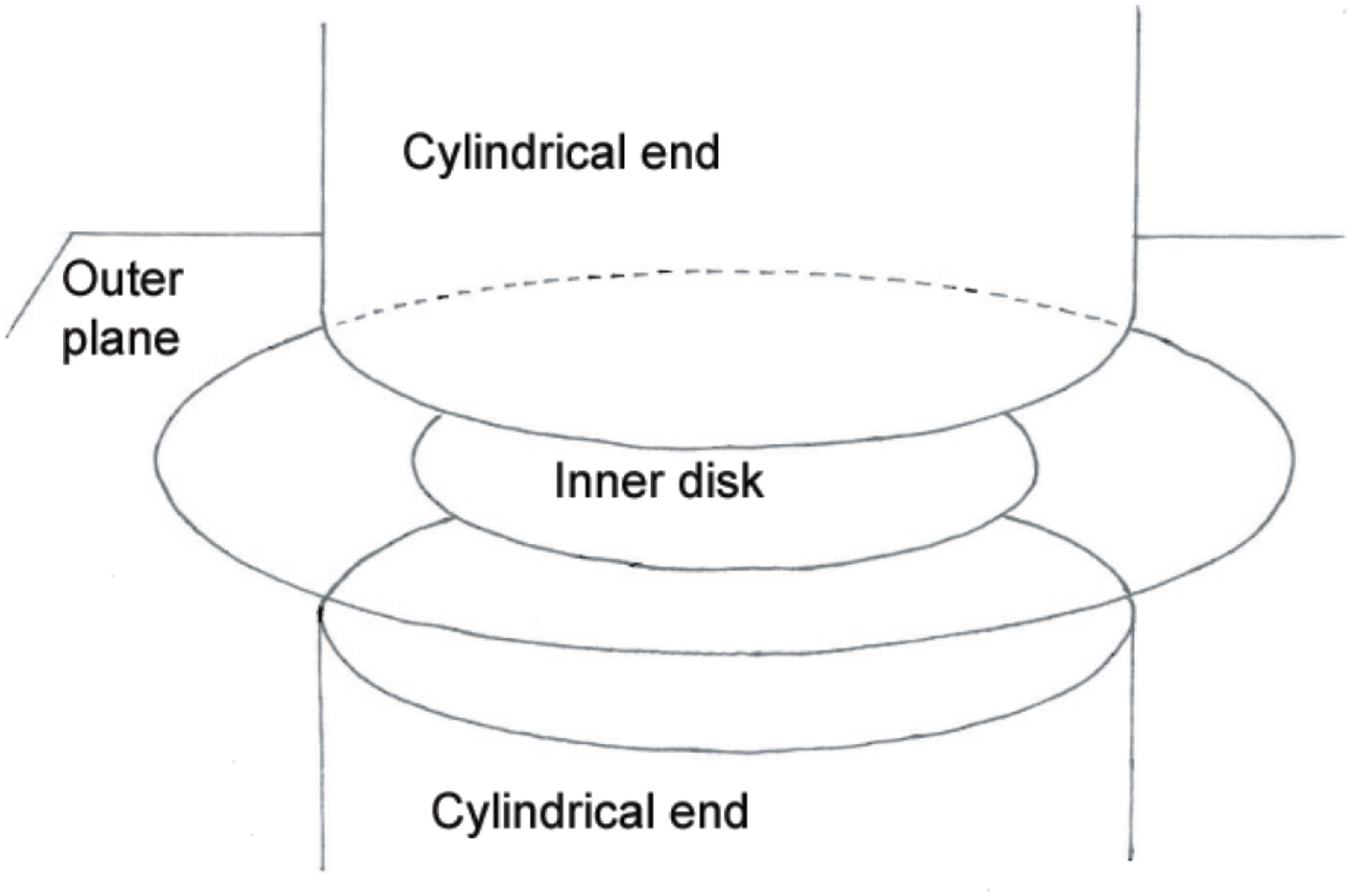}
\hspace*{1cm}
\includegraphics[height=1.4in]{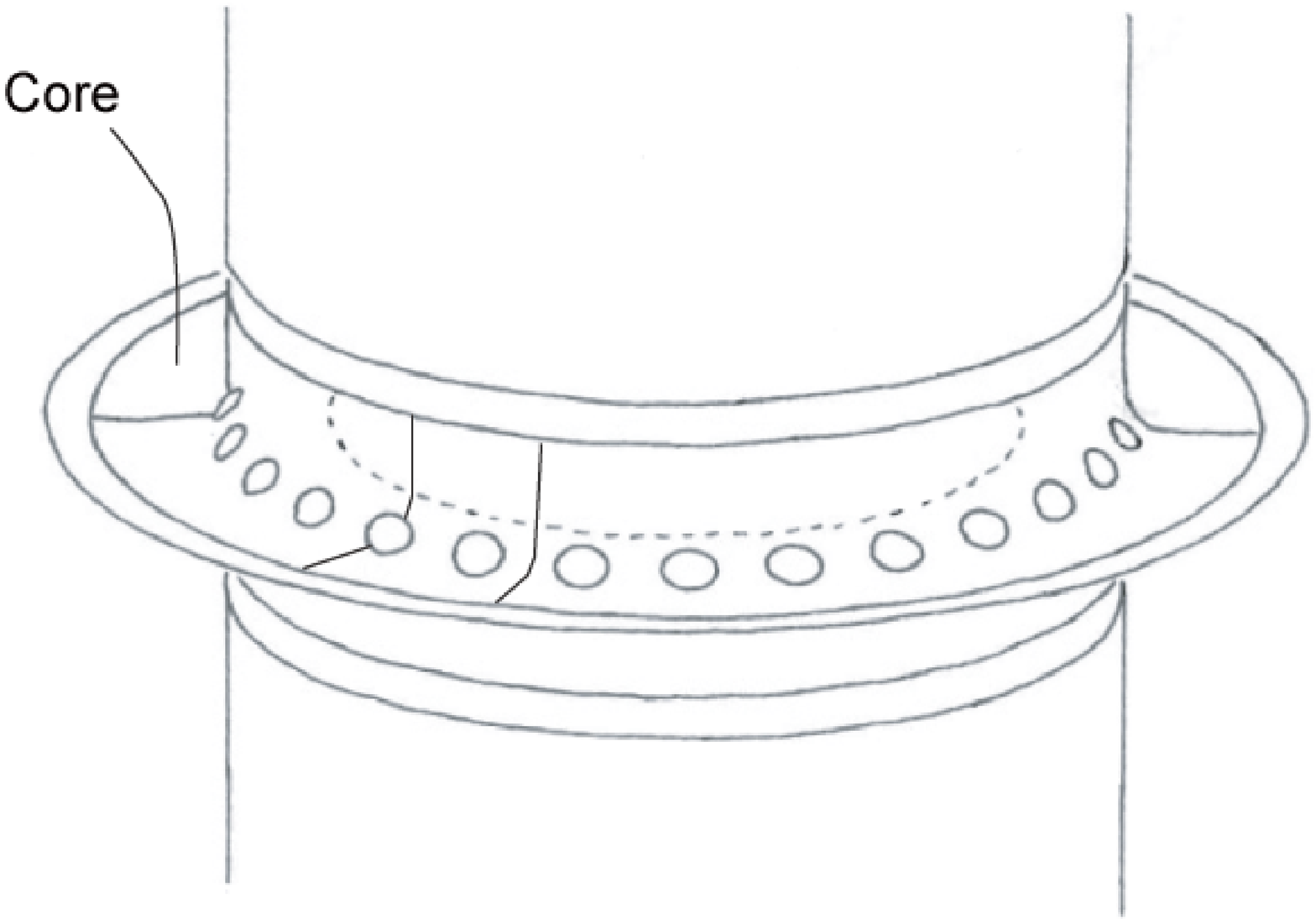}
\caption{Removing a neighborhood of the intersection and inserting the core $\tilde \Sigma^{C}_{1/N}$}
\label{fig:core}
\end{figure}

To construct a new self-similar surface, we take two known examples and replace a neighborhood of their intersection with an appropriately bent scaled Scherk's singly periodic surface, called the core and denoted by $\tilde \Sigma^{C}_{1/N}$.  In the case of a cylinder of radius one and a plane perpendicular to the axis of the cylinder, the desingularization process is  shown in Figure \ref{fig:core}  while  Figure \ref{fig:scherk} shows a portion of the original unmodified Scherk surface. This idea works also if the cylinder is replaced by a sphere.

The resulting surface is not smooth; however it is a good approximate solution. The next task is to find  a function whose graph over it satisfies the self-shrinker equation \eqref{eq:self-shrinker}. Before considering graphs of functions on the entire surface, we have to work locally and study
 the Dirichlet problems with small boundary conditions corresponding to \eqref{eq:self-shrinker}  for graphs of functions on the five different pieces: the core, the outer plane, the inner disk and the two cylindrical ends (or spherical caps, if we work with a sphere). 

In this article, we give a precise description of the maps $\Phi_{\tau}$ used to bend and the maps $\mathcal H_{\tau}$ used to scale the singly periodic Scherk surface as well as tackle the Dirichlet problem on the  central core piece.  The second article shows that the Dirichlet problem corresponding to equation \eqref{eq:self-shrinker} for graphs of functions on the unbounded outer plane has a solution. These are the two more difficult Dirichlet  problems. The results remain valid, whether we use a sphere or a cylinder in the construction. The third article is to discuss the glueing of the solutions on the different pieces in a manner to obtain a smooth complete embedded self-similar surface.

\subsection{Main result} Roughly speaking, our main result states that if the bent scaled Scherk surface used for desingularizing is scaled enough then a small perturbation of it satisfies \eqref{eq:self-shrinker} and  imposed boundary conditions, provided the latter are small enough. The statement is made more precise below.
 
Let us denote the bent scaled surface $\mathcal H_{1/N}( \Phi_{1/N} (S_{C}))$ by $\tilde \Sigma^{C}_{1/N}$ (this is the core in Figure \ref{fig:core}). By our choice of bending and scaling, the surface $\tilde \Sigma^{C}_{1/N}$ has $N$ handles. Moreover, it is invariant with respect to the rotation of $180$ degrees about the $x$-axis and with respect to symmetries across the planes $\theta = \pi/2N + k\pi/N$, $k =1, \ldots, 2N$, where we used cylindrical coordinates.

\begin{theorem}
\label{thm:intro}
There is a constant $C_1$ so that, for any $C>C_1$, there exist an integer  $N_0>0$ and a constant $\delta_0>0$ with the following properties:\\
for every integer $N>N_0$ and for every  function $\tilde f \in W^{2, p}(\tilde \Sigma^C_{1/N})$, with  $\lVert \tilde f \rVert_{W^{2,p}(\tilde \Sigma^C_{1/N})} \leq \delta_0/(2 N^2)$ and whose graph  over $\tilde \Sigma^C_{1/N}$ satisfies the above symmetries, there exists a function $\tilde h \in W^{2,p}(\tilde \Sigma^C_{1/N})$ such that 
	\begin{gather*}
	\text{ the graph of $\tilde h$ over $\tilde \Sigma^C_{1/N}$ satisfies $H+X \cdot \nu =0$,}\\
	\text{ the  graph of $\tilde h$ over $\tilde \Sigma^C_{1/N}$ enjoys the above symmetries,}\\
	 \tilde h=\tilde f \textrm{ on } \pd \tilde \Sigma^C_{1/N}.
	\end{gather*}
Moreover, we can choose $N_0$ and $\delta_0$ so that the graph of $h$ over $\tilde \Sigma^C_{1/N}$ is an embedded surface.
\end{theorem}

 Imposing symmetries on the solutions and boundary conditions that naturally mirror the symmetries of the Scherk surface simplifies the problem greatly. We can  work on a piece of surface characterized by a single period then reconstruct the complete surface using the symmetries. The central tool in the proof is the Inverse Function Theorem. The fact that the Scherk surface is minimal is also used strongly: the Gauss map is then conformal and allows us to transport partial differential equations to the sphere, where elliptic operators have been well studied.

Our methods to solving the Dirichlet problems differ from Kapouleas'. He was able to change slightly the relative position of the Scherk surface wings as well as the size and position of his catenoidal ends in order to compensate for the existence of an approximate kernel. Since in our case, the radius of the cylinder and position of the plane are fixed by  equation \eqref{eq:self-shrinker},  the Dirichlet problems have to be solved directly.

This paper is part of the author's thesis, written under the direction of Sigurd Angenent at the University of Wisconsin-Madison. The author is indebted to professor Angenent for valuable discussion.

\section{Definitions and Bending Maps}
\label{sec:definitions-bendings}

In this section, we will describe in details the transformations made on the Scherk surface $S$ in order to obtain a suitable setting to define and solve a Dirichlet problem. 
\subsection{Scherk surface}   
Consider the Scherk surface $S$ given 
\begin{equation}
\label{eqdef:Scherk-surface}
\sin y = \sinh x \sinh z
\end{equation}
and shown in Figure \ref{fig:scherk}.
\begin{lemma}
\label{lem:Scherk-properties}
The surface $S$ is a singly periodic embedded minimal surface which has the following properties:\\
	- $S$ is asymptotic to the $xy$-plane and to the $yz$-plane,\\
	- $S$ is invariant under rotation of $180^{\circ}$ about the axes $\{ (x,  k \pi, 0) \}$ ($k \in \mathbf{Z}$),\\
	- $S$ is invariant under reflections with respect to the planes $\{y = \pi/2 + k \pi \}$ ($k \in \mathbf{Z}$).
\end{lemma}

\begin{proof}
The symmetries can be readily derived from equation \eqref{eqdef:Scherk-surface}. A direct computation  shows the mean curvature vanishes. For more information on the Scherk  surface, we refer to \cite{dierkes-hildebrandt;minimal-surfaces}.
\end{proof}

We truncate our surface using the  domain 
	\begin{equation}
	\label{eq:truncation}
\Omega=\{(x, y, z) \in \R^3 |\ \  |x| < C, \frac{-3\pi}{2}<y<\frac{\pi}{2}, |z| < C_0\},
	\end{equation}
for some large constant $C_0>100$ chosen later. We omit the dependence on $C_0$ in our notation for $\Omega$ and for the surfaces defined in this article, when we do not wish to emphasize it. 

Fix a constant $\eta$ so that $0<\eta \leq \frac{1}{2C}$.  For all $\tau \in (-2\eta, 2\eta)$, define the bending maps $\Phi_{\tau}:\Omega \to \R^3$ in the following manner:
	\begin{equation}
	\label{eqdef:bendingmap}
	\Phi_{\tau}(x, y, z)=
		\begin{cases}
			\left(\frac{1}{\tau} e^{\tau x} \cos(\tau y) - \frac{1}{\tau}, \frac{1}{\tau} e^{\tau x} \sin(\tau y), z\right), &\quad \tau\neq 0,\\
			(x, y, z), &\quad \tau=0.
		\end{cases} 
	\end{equation}
The maps $\Phi_{\tau}$ transform the $yz-$plane into cylinders of radius $\tau^{-1}$ around the axis $x= -\tau^{-1}, y=0$. Therefore, bending the Scherk surface $S$ using $\Phi_{\tau}$ creates surfaces $\Phi_{\tau}(S)$ asymptotic to cylinders of radius $\tau^{-1}$. In the end, we want our surfaces to be asymptotic to a cylinder of radius $1$; however, working at this larger scale allows us to define a non singular surface for the case $\tau=0$ and a smooth family of surfaces in a neighborhood of $\tau=0$.

Denote by $\Sigma_{\tau}$ the bent truncated Scherk surface
	\begin{equation}
	\label{eqdef:sigmatau}
\Sigma_{\tau}=\Phi_{\tau}(S \cap \Omega).
	\end{equation}
%
%
	\begin{lemma} 
	\label{lem:Phitaudifferentiable}
Let the maps $\Phi_{\tau}$ be  as in (\ref{eqdef:bendingmap}). Then the map 
	\begin{align*}
	\Phi: (-\eta, \eta) & \to C^5(\overline{\Omega}, \R^3)\\
	\tau &\mapsto \Phi(\tau)=\Phi_{\tau}
	\end{align*}
 is continuously Fr\'echet differentiable.
	\end{lemma}

\begin{proof} From the explicit formula \eqref{eqdef:bendingmap} defining the $\Phi_{\tau}$'s, it is immediate that $\Phi$ is continuously differentiable in all its variable up to order six and its derivatives are bounded on $(-\eta, \eta) \times \overline{\Omega}$. The argument is then completed by applying the Mean Value Theorem.
\end{proof}

%

%
\subsection{Tensor bundles on $\Sigma_0$}

Consider $\bar \Sigma_0$, the closure of $\Sigma_0$ as a smooth manifold with boundary. The metric on $\Sigma_0$ is induced by the Euclidean metric on $\R^3$.

A tensor $r$ times contravariant and $s$ times covariant is called an $(r, s)$ type tensor. 
Denote by $T_{r, s}$ the tensor bundle of type $(r,s)$  over $\Sigma_0$. This vector bundle possesses a metric induced by the metric of $\Sigma_0$.  An $(r, s)$ type tensor $f$ over $\Sigma_0$ is a  section of $T_{r, s}$.

\subsection{Definition of $\Sigma_{h, \tau}, X\htau, a\htau, A\htau, \nu\htau$ and $H\htau$}

We consider the surfaces $\overline{\Sigma}_{\tau}$ as smooth manifolds with boundary with respective metrics $g_{\Sigma_\tau}$ induced by the Euclidean metric in $\R^3$. Let $p \in \Sigma_0$ and denote by $X^{(0, 0)}(p)$ and $\nu^{(0, 0)}(p)$ the position vector and the normal vector at $p$ respectively.

The vectors $X^{(0, \tau)}(p)$,  $\nu^{(0, \tau)}(p)$ and the function $H^{(0, \tau)}$ are pull backs of the position vector,  normal vector and the mean curvature at $\Phi_{\tau}(p) \in \Sigma_{\tau}$ to $p \in \Sigma_0$, i.e.
	\begin{align*}
	X^{(0, \tau)}(p) &=  \Phi_{\tau}(X^{(0, 0)}(p)) = \Phi_{\tau}(p)\\
	\nu^{(0, \tau)}(p) & = \textrm{ unit normal to $\Sigma_{\tau}$ at }(\Phi_{\tau}(X^{(0, 0)}(p))\\
	H^{(0, \tau)}(p) & = \textrm{ mean curvature of $\Sigma_{\tau}$ at }(\Phi_{\tau}(X^{(0, 0)}(p)).
	\end{align*}
Similarly, the tensor $g^{(0, \tau)}\in T_{0, 2}$ is the pull back to $\Sigma_0$ of the metric on $\Sigma_{\tau}$, the tensor $a^{(0, \tau)}\in T_{0, 2}$ is the pull back to $\Sigma_0$ of the second fundamental form  on $\Sigma_{\tau}$, and  $A^{(0, \tau)} \in T_{1, 1}$ is the pull back of the second fundamental tensor on $\Sigma_{\tau}$ to $\Sigma_0$ using the map $\Phi_{\tau}$: for $\xi, \rho \in T_p\Sigma_0$,
	\begin{align*}
	g^{(0, \tau)}(p)(\xi, \rho) &= (\Phi_{\tau}^{\ast} g_{\Sigma_{\tau}})(p)(\xi, \rho) = g_{\Sigma_{\tau}}(\Phi_{\tau}(p)) (\Phi_{\tau \ast} \xi, \Phi_{\tau \ast} \rho), \\
	a^{(0, \tau)}(p)(\xi, \rho) &= (\Phi_{\tau}^{\ast} a_{\Sigma_{\tau}})(p) (\xi, \rho) = a_{\Sigma_{\tau}} (\Phi_{\tau}(p))(\Phi_{\tau}^{\ast} \xi, \Phi_{\tau}^{\ast} \rho),\\
	\langle \rho, A^{(0, \tau)}(p)\xi \rangle &= \langle \rho, (\Phi_{\tau}^{\ast} A_{\Sigma_{\tau}})(p)\xi \rangle  = \langle \Phi_{\tau \ast} \rho, A_{\Sigma_{\tau}}(\Phi_{\tau}(p)) (\Phi_{\tau \ast} \xi)\rangle, 	\end{align*}
where $g_{\Sigma_{\tau}}$, $a_{\Sigma_{\tau}}$ and $A_{\Sigma_{\tau}}$ are respectively the metric, the second fundamental form and the type $(1, 1)$ second fundamental tensor of  $(\Sigma_{\tau}, g^{(0, \tau)})$.

For a continuous function $h:\Sigma_0 \to R$, the graph of $h$ over $\Sigma_{\tau}$ is defined by
	\[
	\Sigma_{( h, \tau)} = \{ X^{(0, \tau)}(p) +h(p) \nu^{(0, \tau)}(p), \quad p \in \Sigma_0\}.
	\]

\begin{claim}
\label{claim:embedded}
For all $\eta_0<\frac{1}{200}$, there is a constant $\delta>0$ depending on $\eta_0$ such that for all functions $h:\Sigma_0 \to \R$ with $|h|<\delta$ and constants $\tau \in (-\eta, \eta)$, the surface $\Sigma_{(h, \tau)}$ is embedded.  For example, if $\eta_0 =1/200$, we can take $\delta =1/8$. 
\end{claim}
From the proof below, note that $\delta$ can be chosen larger as the given $\eta_0$ is smaller. \begin{proof}
The transformations $\Phi_{\tau}$ are continuously differentiable with respect to $\tau$ and tend to the identity as $\tau \to 0$. We want the graph of a function on the surface to be embedded; the only problematic parts are the "holes" in the Scherk surface. The "hole" in the domain $0\leq y \leq \pi$ is at its narrowest in the plane $x=z$.  If we denote by $u$ and $y$ the coordinates in this plane, the intersection curve of $S$ with it is characterized by the equation $\sin y = \sinh^2 (u/\sqrt 2)$. The previous equation implies that the graph of a function $h$ with $|h| < 1/2 < \frac{1}{2}\min(\pi, 2 \sqrt 2 \sinh^{-1}1)$ over $\Sigma_0$ is embedded. The mappings $\Phi_{\tau}$ for $\tau \in (-1/100, 1/100)$ only change the surface $\Sigma_0$ slightly and taking $\delta <1/8$ is safe, although it is a very crude estimate.
\end{proof}

For $h$ as in the previous lemma, consider the graphs $\Sigma_{h, \tau}$ as manifolds. We can equip the surfaces $\Sigma_{h, \tau}$ with metrics $g_{\Sigma_{h, \tau}}$ induced by the embedding and the Euclidean metric in $\R^3$. Define 
	\[
	X^{(h, \tau)}(p) = X^{(0, \tau)}(p) +h(p) \nu^{(0, \tau)}(p)
	\]
and the map  $\Phi_{h, \tau}: \Sigma_0 \to \Sigma_{h, \tau}$ by $\Phi_{h, \tau}(p) = X^{(h, \tau)}(p)$. We can now define 
	\begin{align*}
	\nu^{(h, \tau)}(p) & = \textrm{ unit normal to $\Sigma_{h, \tau}$ at } \Phi_{h, \tau}(X^{(0, 0)}(p))\\
	H^{(h, \tau)}(p) & = \textrm{ mean curvature of $\Sigma_{h, \tau}$ at }\Phi_{h, \tau}(X^{(0, 0)}(p)). 
	\end{align*}
The $(0, 2)$ tensor $g\htau$, the $(0, 2)$ tensor $a\htau$ and the $(1,1)$ tensor $A\htau$ are  respectively the metric, the second fundamental form and the second fundamental tensor on $\Sigma_{h, \tau}$  pulled back to $\Sigma_0$ using $\Phi_{h, \tau}$.
	\[
	g\htau = \Phi_{h, \tau}^{\ast} g_{\Sigma_{h, \tau}}, \quad
	a\htau = \Phi_{h, \tau}^{\ast} a_{\Sigma_{h, \tau}}, \quad
	A\htau = \Phi_{h, \tau}^{\ast} A_{\Sigma_{h, \tau}},
	\]
where $g_{\Sigma_{h, \tau}}$, $a_{\Sigma_{h, \tau}}$ and  $A_{\Sigma_{h, \tau}}$ are respectively the metric, the second fundamental form and the second fundamental tensor on $\Sigma_{h, \tau}$.

By the definitions above, all of the vectors and tensors with superscript $(\cdot, \tau)$ are sections of $\Sigma_0$ into either $\R^3$ or a tensor bundle over $\Sigma_0$.

\section{Fr\'echet Differentiability}
\label{sec:F-differentiability}
%
%
\subsection{Scaling, translating and modified equation}

We wish to replace a neighborhood of the  intersection of a cylinder of radius $1$ about the $z$-axis and a plane, by a surface satisfying \eqref{eq:self-shrinker}; however,  the surfaces $\Sigma_{h, \tau}$ are pieces of surfaces asymptotic to cylinders of radius $\tau^{-1}$ about the line $(-\tau^{-1}, 0, z)$. Working in this larger scale requires us to modify the equation above to take into account the different scale and displacement of the axis of the cylinder.

Define the map 
	\begin{equation}
	\label{eqdef:scalingmap}
	\mathcal{H}_{\tau}:\R^3 \to \R^3,\quad \mathcal{H}_{\tau} (x, y, z)= \tau(x+\frac{1}{\tau}, y, z).
	\end{equation}  

\begin{lemma}
\label{lem:scaling-HplusXnu}
If $\tau \neq 0$ and $h \in W^{2,p}(\Sigma_0)$  is such that 
	\begin{equation}
	\label{eq:scaled-eq}
	H\htau+\tau e_1\cdot \nu\htau+\tau^2 X\htau \cdot \nu\htau=0,
	\end{equation}
where $e_1 = (1, 0, 0)$ is the first coordinate vector in $\R^3$, then the rescaled surface $\mathcal{H}_{\tau} (\Sigma_{h, \tau})$ satisfies the equation for contracting self similar surfaces $H + X\cdot \nu =0$,
where $H$, $X$, and $\nu$ are the mean curvature, the position vector and the normal vector to $\mathcal{H}_{\tau}(\Sigma_{h, \tau})$ respectively.
\end{lemma}

\begin{proof}
We have 
	$
	X\htau = \frac{1}{\tau} X - (\frac{1}{\tau}, 0, 0) = \frac{1}{\tau}(X-(1, 0, 0)),
	$
therefore
	$
	H\htau = \tau H,
	\nu\htau = \nu,
	$
and 
	\[
	\begin{split}
	0&=H\htau +\tau e_1\cdot \nu\htau+\tau^2 X\htau \cdot \nu\htau \\
	&= \tau H+ \tau e_1 \cdot \nu +\tau (X - (1, 0, 0))\cdot \nu= \tau (H+X\cdot \nu).\qedhere
	\end{split} 
	\]
	\end{proof}

\subsection{Fr\'echet differentiability}
For $\tau \in (-\eta, \eta)$ and $h \in W^{2,p}(\Sigma_0, \R)$ so that $|h(p)| <1/8$ for $p\in \Sigma_0$, Claim \ref{claim:embedded} guarantees that the graph of $h$ over $\Sigma_{\tau}$ is embedded therefore we can define the function 
	\begin{equation}
	\label{eqdef:F}
	F(h, \tau) = H\htau+\tau e_1\cdot \nu\htau+\tau^2 X\htau \cdot \nu\htau,
	\end{equation}
where $e_1 = (1, 0, 0)$ is the first coordinate vector in $\R^3$. We can now state the main theorem of this section:

\begin{theorem} 
\label{thm:F-differentiable}
Let $p>2$. There is a neighborhood $U \subset W^{2, p}(\Sigma_0)$ of the zero section  so that the map 
	\begin{align*}
	F: U \times (-\eta, \eta) &\to L^p(\Sigma_0)\\
	(h, \tau) &\mapsto H\htau+\tau e_1\cdot \nu\htau+\tau^2 X\htau \cdot \nu\htau
	\end{align*}
 is continuously Fr\'echet differentiable.
\end{theorem}

\subsection{Preliminary lemmas}
Before starting the proof of Theorem \ref{thm:F-differentiable}, we need some preliminary results. The first one shows that the metric, normal vector and second fundamental form of $\Sigma_{\tau}$ are continuously Fr\'echet differentiable ($C^1$) in $\tau$. 
\begin{lemma} 
\label{lem:metric-nu-A}
The map $\tau \mapsto g^{(0, \tau)}$ is $C^1$ from $(-\eta, \eta) \to C^4(T_{0, 2})$.\\
The map $\tau \mapsto g_{(0, \tau)} = (g^{(0, \tau)})^{-1}$ is $C^1$ from $(-\eta, \eta) \to C^4(T_{2, 0})$.\\
The map $\tau \mapsto \nu^{(0, \tau)}$ is $C^1$ from $(-\eta, \eta) \to C^3(\Sigma_0, \R^3)$.\\
The map $\tau \mapsto A^{(0, \tau)}$ is $C^1$ from $(-\eta, \eta) \to C^3(T_{1, 1})$.
\end{lemma}

\begin{proof}
The surface $\Sigma_0$ can be covered by a finite number of coordinate charts,  $\Sigma_0= \bigcup_{l=1}^L O_l$. Using a partition of unity and the fact that  sums and products of $C^1$ functions are  $C^1$, it suffices to show that $\tau \mapsto g^{(0, \tau)}|_{O_l}$ is $C^1$ from $(-\eta, \eta)$ into $C^4(T_{0, 2}|_{O_l})$ for each $l$. Similarly, we only need to show differentiability of the normal vector and the second fundamental form in each $O_l$. These tensors can be written explicitly in coordinates, the result then follows from Lemma \ref{lem:Phitaudifferentiable}.
\end{proof}

The strategy for proving Theorem \ref{thm:F-differentiable} is to write $F$ in terms of tensors that are continuously differentiable in $\tau$ and $h$. We therefore need a few results about differentiability for products, contractions as well as inverse of tensors.

The Sobolev embeddings are true in our context: for $p>2$, $  W^{1,p} (T_{\alpha, \beta})$ $\subset L^{\infty}(T_{\alpha, \beta})$ and there exists a constant $C$ independent of $u$ so that 
	\[
	\| u \|_{L^{\infty} (T_{\alpha, \beta})} \leq C \| u \|_{W^{1, p}(T_{\alpha, \beta})}.
	\]
From this estimate, we have  that 
the product of tensors $(u, v) \to u \otimes v$ is continuously Fr\'echet differentiable
as a map \\
(i)   from $W^{1, p}(T_{\alpha, \beta}) \times W^{1, p}(T_{\gamma, \delta}) \to W^{1, p}(T_{\alpha+\gamma, \beta+\delta})$,\\
(ii) from $L^{p}(T_{\alpha, \beta}) \times L^{\infty}(T_{\gamma, \delta}) \to L^{p}(T_{\alpha+\gamma, \beta+\delta})$.

In addition, a computation in coordinates shows that contraction operations, which convert a tensor of type $(\alpha, \beta)$ to one of type $(\alpha -1, \beta -1)$, are continuously Fr\'echet differentiable from $W^{1, p}(T_{\alpha, \beta}) \to W^{1, p}(T_{\alpha-1, \beta -1})$.

Define $\mathcal{O}_l = \{ u\in W^{1,p}(T_{1, 1}),  |u(p)| < l, \  \forall p\in \Sigma_0\}$. The set $\mathcal{O}_l$ is open in $W^{1, p}(T_{1, 1})$. We have the following lemma for inverse of tensors:

\begin{lemma}
\label{lem:inverse-of-tensors}
The map 
	\begin{eqnarray*} 
	f: \mathcal{O}_{1/2} &\to& W^{1, p}(T_{1, 1})\\
	u &\mapsto & f(u) = (Id - u)^{-1} -Id
	\end{eqnarray*}
is continuously Fr\'echet differentiable.
\end{lemma}

\begin{proof} Note that $ (Id -u)^{-1} = \sum_{k=0}^{\infty}  u^k$ since the series converges for $u \in \mathcal{O}_{1/2}$. This implies that $\lVert(Id -u)^{-1}\rVert_{L^{\infty}} \leq 2$ for $u \in \mathcal{O}_{1/2}$. Moreover, we have $\lVert \nabla f(u)\rVert_{L^p} = \lVert \nabla (Id -u)^{-1} \rVert_{L^p} 	\leq 4 \lVert\nabla u\rVert_{L^p}$. Now, using the formula
\begin{align*}
	\lefteqn{f(u+s) -f(u) = (Id -u-s)^{-1} - (Id -u)^{-1}}\\
		& = (Id-u-s)^{-1} s(Id -u)^{-1} s (Id -u)^{-1} + (Id -u)^{-1} s (Id -u)^{-1},
	\end{align*}
the proof is then straightforward.
\end{proof}

\subsection{Proof of Fr\'echet differentiability}

We are now ready to prove Theorem \ref{thm:F-differentiable}. Recall that the function is  $F(h, \tau) = H\htau + \tau e_1 \cdot \nu\htau + \tau^2 X\htau \cdot \nu \htau$. The mean curvature $H \htau$ is given by $H\htau = g^{-1}_{(h, \tau)} a\htau$ so we have to establish differentiability for the inverse of the metric, the second fundamental form, as well as the normal vector which plays an essential role in the definition of the second fundamental form. \\

\subsubsection{Step 1: Normal vector $\nu\htau$}

Let $p$ be a point in $\Sigma_{\tau}$  and let $\rho$ be a vector tangent to $\Sigma_0$ at $p$. Denote by $D_{\rho}$ the derivative of $\R^3$ valued functions in the direction of $\rho$. A vector  $\tilde \nu \htau = \nu^{(0, \tau)} + \xi$, where $\xi$ is a tangent vector, is normal if 
	\begin{align*}
	0 &= \langle \tilde \nu \htau, D_{\rho} X\htau\rangle = \langle \nu^{(0, \tau)} + \xi, D_{\rho} X^{(0,\tau)} + D_{\rho} h \nu^{(0, \tau)}+ h D_{\rho} \nu^{(0, \tau)}\rangle \\
	&= D_{\rho} h +\langle \xi, (Id - h A^{(0, \tau)}) \rho\rangle= \langle \nabla^{(0, \tau)} h + (Id -h A^{(0, \tau)}) \xi, \rho\rangle,
	\end{align*}
for any tangent vector $\rho$.
The normal direction is therefore given by   $\tilde{\nu}\htau = \nu^{(0, \tau)} -(Id -hA^{(0, \tau)})^{-1} \nabla^{(0, \tau)} h$.

Note that $\nabla^{(0, \tau)} h = \nabla^{(0, 0)} h\  g^{(0, 0)} g^{-1}_{(0, \tau)}$, therefore 
the map $(h, \tau) \mapsto \nabla^{(0, \tau)}  h$ is continuously Fr\'echet differentiable from $W^{2, p}(\Sigma_0, \R)\times (-\eta, \eta) $ to $W^{1, p} (T_{1, 0})$. Using the result on the second fundamental form $A^{(0, \tau)}$ from Lemma \ref{lem:metric-nu-A}, we get that the map
	\[
	\textrm{$(h, \tau) \mapsto hA^{(0, \tau)}$ is $C^1$ from $W^{2, p}(\Sigma_0) \times (-\eta, \eta) \to W^{1, p}(T_{1, 1})$}.
	\]
Denote by $B_{\mu}(0) \subset {W^{2, p}(\Sigma_0)} $  the ball of radius $\mu$  about the zero section. Choose $\mu >0$ so that $ h A^{(0, \tau)} \in \mathcal{O}_{1/2}$ for $h \in B_{\mu}(0)$.  Lemma \ref{lem:inverse-of-tensors} implies that the map
	 $(h, \tau) \mapsto ((Id-hA^{(0, \tau)})^{-1} -Id)\nabla h$ is $C^1$ from $B_{\mu}(0) \times (-\eta, \eta)\to W^{1, p}(T_{1, 0})$.
Hence, 
	\[
	\textrm{$(h, \tau) \mapsto \tilde{\nu}\htau - \nu{(0, \tau)}$ is $C^1$ from $B_{\mu}(0)\times (-\eta, \eta) \to W^{1, p}(T_{1, 0})$.}
	\]
In order to have the normal vector $\nu\htau$ instead of the normal direction $\tilde{\nu}\htau$, define the function $\tilde b:= |\tilde{\nu}\htau|^{-1}$ on $\Sigma_0$ and write
	\begin{equation}
	\label{eq:nuhtau}
	\nu\htau -\nu^{(0, \tau)} =\tilde b\ \tilde{\nu}\htau -\nu^{(0, \tau)}= \tilde b\ (\tilde{\nu}\htau -\nu^{(0, \tau)}) -(\tilde b-1) \nu^{(0, \tau)}. 
	\end{equation}
Since $\tilde{\nu}\htau -\nu^{(0, \tau)}$ is a tangent vector, 
	$
	\tilde b^{-2}
	=  |\tilde{\nu}\htau-\nu^{(0, \tau)}|^2+1$.
Hence
	\[ 
	\tilde b-1 = \frac{1}{\sqrt{ |\tilde{\nu}\htau-\nu^{(0, \tau)}|^2+1}} -1.
	\]
Thanks to the differentiability of $\tilde{\nu}\htau-\nu^{(0, \tau)}$ and the fact that products and contractions of tensors are differentiable in $W^{1,p}$,  the map $|\tilde{\nu}\htau-\nu^{(0, \tau)}|^2$ is continuously Fr\'echet differentiable from $B_{\mu}(0) \times (-\eta, \eta)$ to $W^{1, p}(\Sigma_0)$ or $L^{\infty}(\Sigma_0)$. Since, $\tilde \nu^{(0, \tau)} = \nu^{(0, \tau)}$, choosing $\mu>0$ smaller if necessary, we can assume without loss of generality that $|\tilde{\nu}\htau-\nu^{(0, \tau)}|^2 < 1/2$ for $h \in B_{\mu}(0)$ and $\tau \in (-\eta, \eta)$.
The following lemma shows that $\tilde b-1$ is differentiable.

\begin{lemma}
\label{lem:sub}
Let $s$ be a section of $T_{0, 0}$, in other words $s$ is a function from $\Sigma$ to $\R$, such that $|s(p)| \leq 1/2 $ for $p\in \Sigma_0$. The map $	W^{1, p}(T_{0, 0})  \to  W^{1, p}(T_{0, 0})$, $
	s \mapsto f \circ s   := \frac{1}{\sqrt{s+1}} -1
	$
is continuously Fr\'echet differentiable.
\end{lemma}
\begin{proof}
Simple estimates yield that $f\circ s \in W^{1, p}(T_{0, 0})$.
Noting that the function $f: (-1/2, 1/2) \to \R$ is smooth and has uniformly bounded derivatives, we use the Mean Value Theorem to conclude. 
\end{proof}

Lemma \ref{lem:sub} implies that $\tilde b-1$ is continuously Fr\'echet differentiable from $B_{\mu}(0) \times (-\eta, \eta) \to W^{1, p}(\Sigma_0)$. Since $\Sigma_0$ is bounded, $\tilde b$ is also $C^1$. Using equation \eqref{eq:nuhtau} and differentiability of products of tensors from $W^{1,p} \times W^{1,p}$ to $W^{1,p}$, we have that the map
	\begin{equation}
	\label{eq:nu-differentiable}
	\textrm{$(h, \tau) \mapsto \nu \htau$ is $C^1$ from $B_{\mu} (0) \times (-\eta, \eta) \to W^{1, p}(\Sigma_0)$}.
	\end{equation}
\subsubsection{Step 2: Inverse of metric}

Let $p$ be a point in $\Sigma_{\tau}$  and let $\xi$ be a vector tangent to $\Sigma_0$ at $p$. Denote by $D_{\xi}$ the derivative  in the direction of $\xi$. Then, 
	\begin{align*}
	D_{\xi} X\htau = D_{\xi} (X^{(0, \tau)} + h \nu^{(0, \tau)}) 
			 = (Id -h A^{(0, \tau)}) \xi + D_{\xi} h\ \nu^{(0, \tau)}
	\end{align*}
and
	\begin{multline*}
	g\htau (\rho, \xi)  =  \langle D_{\rho} X\htau, D_{\xi} X\htau\rangle\\
		= g^{(0, \tau)} (\rho, \xi) - 2 h \langle \rho, A^{(0, \tau)}\xi\rangle  +  h^2 \langle A^{(0, \tau)}\rho, A^{(0, \tau)}\xi \rangle  + D_{\rho} h D_{\xi} h.
	\end{multline*}
	
Define the tensors $a(\rho, \xi) = \langle \rho, A^{(0, \tau)} \xi\rangle$, $b(\rho, \xi) = \langle A^{(0, \tau)}\rho, A^{(0, \tau)} \xi\rangle$ on $\Sigma_0$. Since $a$ and $b$ do not depend on $h$, Lemma \ref{lem:metric-nu-A} implies that $ha$ and $h^2 b$ are continuously Fr\'echet differentiable from $B_{\mu}\times (-\eta, \eta) \to C^1(T_{0, 2})$.
Therefore, the map $(h, \tau) \mapsto g\htau - g^{(0, \tau)}$ is $C^1$ from $B_{\mu}(0)\times (-\eta, \eta) \to W^{1, p} (T_{0, 2})$. Let
	$
	u= g^{-1}_{(0, \tau)} (g^{(0, \tau)}-g\htau ).
	$
Choosing $\mu$ smaller if necessary, we can apply Lemma \ref{lem:inverse-of-tensors} to $u$ and get that the map
	\[
	(Id -u)^{-1} -Id = [Id - g^{-1}_{(0, \tau)} (g^{(0, \tau)}-g\htau )]^{-1} - Id = (g^{-1}_{(h, \tau)}-g^{-1}_{(0, \tau)}  ) g^{(0, \tau)}
	\]
is $C^1$ as a map from $B_{\mu}(0)\times (-\eta, \eta)$ to $W^{1, p}(T_{1, 1})$. Finally, we have that the map $(h, \tau) \mapsto g^{-1}_{(h, \tau)}-g^{-1}_{(0, \tau)} $ is $C^1$ from $B_{\mu}(0)\times (-\eta, \eta) \to W^{1, p}(T_{2, 0})$,
therefore the map 
	\begin{equation}
	\label{eq:metric-differentiable}
	\textrm{
	$(h, \tau) \mapsto g^{-1}_{(h, \tau)} $ is $C^1$ from $B_{\mu}(0)\times (-\eta, \eta)  \to W^{1, p}(T_{2, 0})$.}
	\end{equation}

\subsubsection{Step 3: Second fundamental form}

We start with a computation giving us an equation for the second fundamental form $a\htau$ in terms of known quantities. Recall $X\htau = X^{(0, \tau)} +h \nu^{(0, \tau)}$. Denote by $ \nabla^{(0, \tau)}$ the covariant derivative and by $(\nabla^{(0, \tau)})^2$ the second covariant derivative on $\Sigma_0$ corresponding to the metric $g^{(0, \tau)}$. For $\xi, \rho$ tangent vectors to $\Sigma_0$ at the point $p$, we have
	\begin{multline*}
	(\nabla^{(0,  \tau)})^2 X\htau (\rho, \xi) =
	 (\nabla^{(0,  \tau)})^2 X^{(0, \tau)}(\rho, \xi) + h (\nabla^{ (0, \tau)})^2\nu^{(0, \tau)} (\rho, \xi) \\- A^{(0, \tau)}\rho D_{\xi} h - A^{(0, \tau)}\xi D_{\rho} h +(\nabla^{(0, \tau)})^2 h (\rho, \xi)\nu^{(0, \tau)}.
	\end{multline*}
The second derivative of $\nu^{(0, \tau)}$ is given by 
	\begin{equation}
	\label{eq:second-deriv-nu}
	(\nabla^{ (0, \tau)})^2 \nu^{(0, \tau)} (\rho, \xi) 
				= -\nabla^{ (0, \tau)}_{\rho} (A^{(0, \tau)}) \xi- [\langle A^{(0, \tau)} \rho, A^{(0, \tau)}\xi \rangle] \nu^{(0, \tau)}.
	\end{equation}
The second fundamental form $a^{(h, \tau)}$ is given by $a\htau = |\tilde{\nu}\htau|^{-1} \tilde{a}^{(h, \tau)}$, where  $\tilde{\nu}\htau = \nu^{(0, \tau)} -(Id -hA)^{-1} \nabla^{(0, \tau)} h$ is the normal direction on $\Sigma_{h, \tau}$ and 
	\begin{align*}
	 \tilde{a}^{(h, \tau)} (\rho, \xi) &= |\tilde{\nu}\htau|\langle \rho, A^{(h, \tau)} \xi\rangle \\
	 & = \langle \tilde{\nu}\htau, (\nabla^{(h, \tau)})^2 X\htau (\rho, \xi) \rangle\\
	 & = \langle \tilde{\nu}\htau, (\nabla^{(0, \tau)})^2 X\htau (\rho, \xi) \rangle.
	\end{align*}
The last equality is valid since $D_{\xi} X\htau$ is perpendicular to $\tilde{\nu}\htau$. Hence,
	\begin{multline}
	\label{eq:tilde-a}
	\tilde{a}^{(h, \tau)} (\rho, \xi) 
		= a^{(0, \tau)} (\rho, \xi) - h \langle A^{(0, \tau)}\rho, A^{(0, \tau)}\xi\rangle  + (\nabla^{(0, \tau)})^2 h (\rho, \xi) \\
		 + \langle (Id - h A^{(0, \tau)})^{-1} \nabla^{(0, \tau)} h,  h \nabla^{(0,\tau)}_{\rho} A^{(0, \tau)} \xi+  A^{(0, \tau)}\rho D_{\xi} h+ A^{(0, \tau)}{\xi} D_{\rho} h\rangle.
	\end{multline}
We know that 
	the map $(h, \tau) \mapsto (\nabla^{(0, 0)})^2 h$ is $C^1$ from $B_{\mu}(0) \times (-\eta, \eta) \to L^p(T_{0, 2})$.
Let us denote by $(\Gamma^{(0, \tau)})^k_{ij}$ and $(\Gamma^{(0, 0)})^k_{ij}$ the components of the Christoffel symbols corresponding to the metrics $g^{(0, \tau)}$ and $g^{(0, 0)}$ respectively. By Lemma \ref{lem:metric-nu-A}, the $(1, 2)$ tensor given in coordinates by $(\Gamma^{(0, \tau)})^k_{ij}-(\Gamma^{(0, 0)})^k_{ij}$ is $C^1$ from $B_{\mu}(0) \times (-\eta, \eta) \to C^3(T_{1, 2})$. Therefore
	\begin{align*}
	\textrm{
	the map }(h, \tau) &\mapsto (\nabla^{(0, \tau)})^2 h \textrm{ is $C^1$ from $B_{\mu}(0) \times (-\eta, \eta) \to L^p(T_{0, 2})$,}\\
	\textrm{
	the map }(h, \tau) &\mapsto \nabla^{(0, \tau)} A^{(0, \tau)} \textrm{ is $C^1$ from $B_{\mu}(0) \times (-\eta, \eta) \to C^2(T_{1, 2})$.}
	\end{align*}
The other quantities in equation \eqref{eq:tilde-a} have been discussed previously, so it follows from the results on products and contractions of tensors  that 
	the map $(h, \tau) \mapsto  \tilde a^{(h, \tau)}$ is $C^1$ from $B_{\mu}(0) \times (-\eta, \eta) \to L^p(T_{0, 2})$. Since $ |\tilde{\nu}\htau|^{-1} $ is a continuously Fr\'echet differentiable function of $(h, \tau)\in B_{\mu}(0) \times (-\eta, \eta)$ into $W^{1, p}(\Sigma_0)$, 
	\begin{equation}
	\label{eq:a-differentiable}
	\textrm{
	the map }(h, \tau) \mapsto a^{(h, \tau)} \textrm{ is $C^1$ from $B_{\mu}(0) \times (-\eta, \eta) \to L^p(T_{0, 2})$.}
	\end{equation}

\subsubsection{Final step}

Combining the results for the normal vector, the metric and the second fundamental form \eqref{eq:nu-differentiable} \eqref{eq:metric-differentiable} and \eqref{eq:a-differentiable}, we have that
	\[
	\textrm{
	the map }(h, \tau) \mapsto F(h, \tau) \textrm{ is $C^1$ from $B_{\mu}(0) \times (-\eta, \eta) \to L^p(T_{0, 2})$.}
	\]

\section{Inverse Function Theorem Argument}
\label{sec:inverse-f-argument}

Our main goal in this section is to  prove that, loosely speaking, for small enough boundary conditions on $\pd \Sigma_0 \cap (\{ x=\pm C_0\} \cup \{ z= \pm C_0\})$ and small $\tau$, there is a function $h$ with the boundary conditions and so that $F(h, \tau) =0$. We will impose symmetries on the boundary conditions and the function $h$ in order to simplify the problem and to be able to reconstruct a complete surface from the piece $\Sigma_{h, \tau}$. The symmetries do not allow us to choose the boundary values on the curves $y=\pi/2 \cap \pd \Sigma_0$  or $y = 3\pi/2 \cap \pd \Sigma_0$.  The following theorem states the result more precisely.
\begin{theorem}
\label{thm:existence-of-h}
There exists a constant $C_1$ so that, for any $C_0>C_1$, there are constants $\eta_0>0$ and $\delta_0>0$ such that, if we denote by  $\Sigma_0$ the surface $
	\Sigma_0 = \{ (x, y, z) | \sin y =\sinh x \sinh z, -C_0 < x, z < C_0, -\pi/2 <y<3\pi/2\},
	$
then for any $\tau\in (-\eta_0, \eta_0)$ and  function $f \in W^{2, p}(\Sigma_0)$ satisfying $\lVert f \rVert_{W^{2,p}(\Sigma_0)} \leq \delta_0$ and $f(x, y, z) =-f(x, -y, -z) = f(x, \pi -y, z)$, $(x, y, z) \in \Sigma_0$,
there exists a function $h \in W^{2, p}(\Sigma_0)$ having the same symmetries as $f$ and satisfying
	\begin{gather*}
	F(h, \tau) = H\htau+\tau e_1\cdot \nu\htau+\tau^2 X\htau \cdot \nu\htau= 0, \\
	 h=f \textrm{ on } \pd \Sigma_0 \cap (\{ x=\pm C_0\} \cup \{ z= \pm C_0\}).
	\end{gather*}
Moreover, $\eta_0$ and $\delta_0$ can be chosen small enough so that the graph of $h$ over $\Sigma_{\tau}$, $\tau\in (-\eta_0, \eta_0)$ is embedded.
\end{theorem}
The proof is given in Section \ref{ssec:corollaries}. This theorem is a direct consequence of Theorem \ref{thm:G-diffeomorphism} below, which uses an inverse function theorem argument. Before we can state it, we need to define the following Banach spaces. 

\subsection{$W^{2, p}_{Sym}(\Sigma_0)$, Sobolev space with symmetries}
\label{ssec:sobolev-sym}

To simplify our problem, we impose symmetries on the functions $h$. Roughly speaking, we will only consider functions $h$ such that the graph $\Sigma_{h, 0}$ of $h$ over $\Sigma_0$ satisfies the same symmetries as $\Sigma_0$. More precisely, $h \in W^{2, p}_{Sym}$ if $\Sigma_{h, 0}$ is invariant under \\
	- rotation of $180^{\circ}$ with respect to the  $x$-axis: $h(x, -y, -z) = - h(x, y, z)$, $(x, y, z) \in \Sigma_0$,\\
	- symmetry with respect to the plane $y = \pi/2$: $h(x, y, z) = h(x, \pi-y, z)$, $(x, y, z) \in \Sigma_0 $.\\
Let us keep in mind that the two symmetries above, considered on the whole Scherk surface $S$, generate all the symmetries listed in Lemma \ref{lem:Scherk-properties}. Define the Sobolev space with symmetries
	\[
	 W^{2, p}_{Sym}(\Sigma_0) = \{ h \in W^{2, p}(\Sigma_0) |\  h \textrm{ satisfies the symmetries above}\}		\]
$W^{2, p}_{Sym}(\Sigma_0)$ ($p>2$) is a Banach space. Since $W^{2, p}_{Sym}(\Sigma_0) \subset C^1(\Sigma_0)$ it makes sense to consider derivatives of $h$. 

\subsection{ Definition of $\Sigma_0^{\ast}$} 

Let us identify pairs of points on the edges $\pd \Sigma_0 \cap\{y=-\pi/2, 3\pi/2\}$ and denote by $\Sigma_0^{\ast}$ the set $\Sigma_0$, where  the points $(x, -\pi/2, z)$ and  $(x, 3\pi/ 2, z)$ are identified.  The metric on $\Sigma_0^{\ast}$ is induced by the metric on $\Sigma_0$.

The second  symmetry condition in the previous section, $h(x, y, z) = h(x, \pi-y, z)$, allows us to make sense of a function $h \in W^{2,p}_{Sym}$ as a function on $\Sigma_0^{\ast}$. It is clear that  
 	\[
	W^{2, p}_{Sym}(\Sigma_0)\cong W^{2, p}_{Sym}(\Sigma_0^{\ast}) =: W^{2,p}_{Sym}.
	\]
The next section shows that $\Sigma^{\ast}$ can be conformally mapped to a part of the $2$- sphere. 

\subsection{Going from $\Sigma_0^{\ast}$ to the sphere}
\label{sec:Sigma-to-sphere}

Let $S$ be the Scherk surface defined by \eqref{eqdef:Scherk-surface}. The following results about the Gauss map of the Scherk surface $S$ are standard (see \cite{kapouleas;embedded-minimal-surfaces} or \cite{dierkes-hildebrandt;minimal-surfaces}).

\begin{lemma}
\label{lem:Gauss-map}
The Gauss map $\nu$ of $S$ has the following properties:\\
(i) $\nu$ restricts to a diffeomorphism from $S \cap \{ y \in [-\pi/2, \pi/2]\}$ onto $S^2 \cap \{  y \geq 0\} \setminus \{ (\pm 1, 0, 0), (0, 0, \pm 1)\}$.\\
(ii) Let $E_i (i=1, \ldots, 4)$ be the arcs into which the equator $S^2 \cap \{ y =0\}$ is decomposed by removing the points $(\pm1, 0, 0)$ and $(0, 0, \pm1)$, numbered so that $E_i$ is in the $i$-th quadrant of the $xz$-plane.
	We then have
		\[
		\nu(S \cap \{ y =-\pi/2\}) = E_2 \cup E_4, \quad  \nu(S \cap \{ y =\pi/2\}) = E_1 \cup E_3.
		\]
(iii) $S$ has no umbilics and $\nu^{\ast} g_{S^2} = \frac{1}{2} |A|^2 g.$
\end{lemma}
$\Sigma_0$ is a piece of a Scherk minimal surface. Since $\Sigma _0$ does not have any umbilics, the Gauss map $\nu: \Sigma_0 \to S^2$, $p \mapsto \nu(p)$ is conformal.  Moreover, we can compute the normal vector $\nu$ explicitly using  equation (\ref{eqdef:Scherk-surface}) and find
	\[
	\nu(x, y, z) = \left( \frac{\sinh z}{\cosh z}, -\frac{-\cos y}{\cosh z \cosh x}, \frac{\sinh x}{\cosh x}\right), \textrm{ for } (x, y, z) \in \Sigma_0.
	\]
The images under the Gauss map of the curves $\{z = \pm C_0\} \cap \pd \Sigma_0$, $\{x = \pm C_0\} \cap \pd \Sigma_0$ are circles on $S^2$ centered at $(\pm 1, 0, 0)$ and $(0, 0, \pm 1)$ respectively. 

Denote by $S^2_{\varphi_0}$ the sphere $S^2$ minus four balls of radius $\varphi_0$ centered at $(\pm 1, 0, 0)$ and $(0, 0, \pm 1)$:
	\[
	S^2_{\varphi_0} = \{ (x, y, z) \in S^2 | x, z \in (-\cos \varphi_0, \cos \varphi_0)\}.
	\]
	The symmetries and identification made on $\Sigma_0$ to obtain $\Sigma_0^{\ast}$  allow us to define a bijective conformal continuous Gauss map $\nu$ from $\Sigma_0^{\ast} $ to $S^2_{\varphi_0}$, with $\varphi_0= \cos^{-1}(\frac{\sinh C_0}{\cosh C_0})$. 
	
	 The rotation about the $x$-axis in $\Sigma_0$ corresponds to the reflection across the $yz$-plane in $S^2 \subset \R^3$, and the reflections across the plane $y= \pi/2$ of $\Sigma_0$ corresponds to reflection across the $xz$-plane in $S^2$. Therefore, if we define
	 $
	 W^{2, p}_{Sym}(S^2_{\varphi_0}) =\{ h' \in W^{2, p}(S^2_{\varphi_0}) | h'(x, y, z) = - h'(-x, y, z) = h'(x, -y, z)\}$, 
where the $W^{2, p}$ norm is taken with respect to the standard metric $g_{S^2}$ on $S^2_{\varphi_0}$, then
	\[
	W^{2, p}_{Sym}(\Sigma_0^{\ast}) \cong W^{2, p}_{Sym}(S^2_{\varphi_0}).
	\]

\subsection{$\mathcal{T}$, the space of  traces of $W^{2, p}_{Sym}(\Sigma_0)$}

 Denote by $\Sigma_0^{\ast}$ the set $\Sigma_0$ where we have identified the points $(x, -\pi/2, z)$ and  $(x, 3\pi/ 2, z)$. We say that 
 	\[
	f, h \in W^{2, p}_{Sym} \textrm{ have the same trace if } f-h \in W^{1, p}_0(\Sigma_0^{\ast})\cap W^{2, p}_{Sym},
	\]
where $W^{1, p}_0(\Sigma_0^{\ast})$ is the closure in the $W^{1, p}$ norm of smooth compactly supported functions on $\Sigma_0^{\ast}$. Define the Banach space $\mathcal{T}$  of traces of $W^{2,p}_{Sym}$ and its norm by
	\begin{gather*}
	\mathcal{T} \cong W^{2, p}_{Sym} /(W^{1,p}_0 (\Sigma_0^{\ast})\cap W^{2, p}_{Sym})\\
	\lVert [f]\rVert_{\mathcal{T}} = \min_{u-f \in W^{1,p}_0 (\Sigma_0^{\ast})\cap W^{2, p}_{Sym}} \lVert u\rVert_{W^{2, p}},\quad  [f] \in \mathcal{T}.
	\end{gather*}
The minimum is achieved since $W^{2,p}$ is a reflexive Banach space. Note that the trace characterizes the value of a function on $x, z =\pm C_0$. The parts of the boundary of $\Sigma_0$ on $y = -\pi/2, 3 \pi /2$ are subject to our symmetry conditions,  we therefore do not have the freedom to impose boundary conditions on them. 
\begin{lemma}
\label{lem:Dirichlet-problem-Laplacian}
For $[f] \in \mathcal{T}$, there is a unique $u \in W^{2, p}_{Sym} $ so that 
	\begin{equation}
	\label{eq:Dirichlet-problem-Laplacian}
	\begin{array}{ll}
	\Delta u =0 \textrm{ in } \Sigma_0^{\ast}\\
	u-f \in W^{1, p}_0(\Sigma_0^{\ast}),
	\end{array}
	\end{equation}
where $\Delta$ is the Laplacian computed using the metric on $\Sigma_0^{\ast}$ induced by $g^{(0, 0)}$ on $\Sigma_0$. 
\end{lemma}
\begin{proof} 
Using the continuous conformal Gauss map $\nu : \Sigma_0^{\ast} \to S^2_{\varphi_0}$, we can consider the following equivalent problem on the sphere:
	\begin{eqnarray*}
	\Delta_{g_{S^2}} v =0 \textrm{ in } S^2_{\varphi_0}\\
	v-f \circ \nu^{-1} \in W^{1, p}_0(S^2_{\varphi_0}),
	\end{eqnarray*}
which is known to have a unique solution $v \in W^{2, p}_{Sym}(S^2_{\varphi_0})$. The pull-back $\nu^{\ast} g_{S^2}$ of the metric on $S^2$ is $\nu^{\ast} g_{S^2} = \frac{1}{2} |A|^2 g$ so function
	$
	u = \nu^{\ast} v = v \circ \nu
	$
is a solution to (\ref{eq:Dirichlet-problem-Laplacian}). It is unique since $v$ is unique.	 
\end{proof}

\subsection{Statement of theorem}
Given $C_0$ a large constant, denote by $\Sigma_0$ the surface 
	$
	\Sigma_0 = \{ (x, y, z) | \sin y =\sinh x \sinh z, -C_0 < x, z < C_0, -\pi/2 <y<3\pi/2\}.
	$
Recall that $F(h, \tau) =H\htau+\tau e_1\cdot \nu\htau+\tau^2 X\htau \cdot \nu\htau$. Let $W^{2,p}_{Sym}(\Sigma_0)$ and $L^p_{Sym}$ be the Sobolev spaces with symmetries defined in section \ref{ssec:sobolev-sym} and  
let $\mathcal{T}$  be the space of traces of $W^{2, p}_{Sym}(\Sigma_0)$ functions. 
	
\begin{theorem}
\label{thm:G-diffeomorphism}
There is a constant $C_1$ such that for all $C_0>C_1$, there exist a constant $\eta_0>0$, a neighborhood $U_{Sym} \subset W^{2, p}_{Sym}(\Sigma_0)$ of the zero section and a neighborhood $V \subset   L^p_{Sym} \times \mathcal{T} \times (-\eta_0, \eta_0)$ of $(0, [0], 0)$ such that the map
	\begin{eqnarray*}
	G: U_{Sym} \times (-\eta_0, \eta_0)& \to &  L^p_{Sym} \times \mathcal{T} \times (-\eta_0, \eta_0) \\
	(h, \tau) & \mapsto & (F(h, \tau), [h], \tau), 
	\end{eqnarray*}
is a diffeomorphism from $U_{Sym} \times (-\eta_0, \eta_0)$ to $V$.
\end{theorem}

Note that $G(0, 0) = (0, [0], 0)$. The proof uses the Inverse Function Theorem.  The map $(h, \tau) \mapsto [h]$ is continuous and linear, therefore continuously Fr\'echet differentiable. Since $F(h, \tau)$ is $C^1$ in its variables, $G$ is continuously differentiable also. In order to invoke the Inverse Function Theorem, we have to show that the derivative of $G$ at $(0, 0)$ is an isomorphism. We compute the elliptic operator characterizing $DG(0, 0)$ below. In section \ref{ssec:eigen}, we use the conformal Gauss map to transport the operator to the $2-$sphere where results on the Laplacian are readily available. The image of $\Sigma_0^{\ast}$ under the Gauss map is $S^2_{\varphi_0}$, a sphere minus four small balls. Loosely speaking, since this surface is close to $S^2$,  the eigenvalues of the Laplacian on $S^2_{\varphi_0}$ are close to its eigenvalues on $S^2$. Finally, we show that the operator on $S^2_{\varphi_0}$ which is equivalent to $DG(0, 0)$  is invertible.

Let us compute $DG(0, 0)$. 
	\begin{align*}
	DG(0, 0) (u, \sigma) &= DG(0, 0) (u, 0) + DG(0, 0)(0, \sigma)\\
		&= (\Delta_{g^{(0, 0)}} u +|A^{(0, 0)}|^2 u, [u], 0) + DG(0, 0) (0, \sigma).
	\end{align*}
Recall that 
	$
	F(0, \tau) = H^{(0, \tau)} + \tau e_1 \cdot \nu^{(0, \tau)} + \tau^2 X^{(0, \tau)} \cdot \nu^{(0, \tau)},
	$
therefore, 
	\begin{align*}
	DF(0, 0) (0, \sigma)  &= \frac{d}{dt} H(0, t\sigma) |_{t=0} + \sigma e_1 \cdot \nu^{(0, 0)}\\
		 & = D_{\tau} H(0, 0) \sigma +  \sigma e_1 \cdot \nu^{(0, 0)}\\
		& = \sigma v, 
	\end{align*}
where $v \in L^p_{Sym}(\Sigma_0^{\ast})$ (in fact, $v \in C^1(\Sigma_0)$ by Lemma \ref{lem:metric-nu-A}). Hence,
	\[
	DG(0, 0) (u, \sigma) = (\Delta u + |A|^2 u + \sigma v, [u], \sigma),
	\]
for some function $v \in L^p_{Sym}$ independent of $(u, \sigma)$. We  discuss the properties of the operator $\Delta + |A|^2$ in the next section.

\subsection{Eigenvalues and eigenfunctions of the operator $\Delta +|A|^2$}
\label{ssec:eigen}

Recall that the Gauss map $\nu$ is conformal from $\Sigma_0^{\ast} \to S^2_{\varphi_0}$ and 
	$
	\nu^{\ast} g_{S^2} = \frac{1}{2} |A|^2 g.
	$
Therefore, 
	\[
	\Delta_{\nu^{\ast} g_{S^2}} +2 = \frac{2}{|A|^2} \Delta_g +2 = \frac{2}{|A|^2} (\Delta_g + |A|^2).
	\]
 To show that the kernel of $\Delta_g + |A|^2$ is trivial in $W^{2, p}_{Sym}(\Sigma_0^{\ast} )\cap W^{1, p}_0 (\Sigma_0^{\ast})$, it suffices to show that problem 
	\begin{eqnarray*}
	\Delta_{g_{S^2}} \zeta + 2 \zeta =0 \textrm { in } S^2_{\varphi_0}\\
	\zeta \in W^{1, p}_0(S^2_{\varphi_0}) \cap W^{2, p}_{Sym}(S^2_{\varphi_0})
	\end{eqnarray*}
only admits the zero function as solution.

Eigenvalues and eigenfunctions of the Laplace operator on the whole sphere have been well studied and are  known (see \cite{montiel-ros;operators-holomorphic}, for example). Denote by $\mathcal R_{xz} $ the reflection with respect to $xz$-plane and by $\mathcal R_{yz}$ the reflection with respect to the $yz$-plane. We are interested in functions on the sphere that are $\mathcal R_{xz}$ invariant and $\mathcal R_{yz}$ antivariant (see discussion from section \ref{sec:Sigma-to-sphere}).   The first  eigenvalues and eigenfunctions of $M_1 = S^2$ are listed in the table below. \\

	\begin{tabular}{c|c|c}
	Eigenvalue & Eigenfunction & Properties of eigenfunction\\ \hline
	$\lambda_{1, 1} =0$ &$ u_{1, 1}$ & $\mathcal R_{xz}$ invariant \\
	$\lambda_{1, 2} =2$ & $u_{1, 2} = y$ & $\mathcal R_{xz}$ antivariant\\
	$\lambda_{1, 3} =2$ & $u_{1, 3} = x$ & $\mathcal R_{xz}$ invariant, $\mathcal R_{yz}$ antivariant\\
	$\lambda_{1, 4} =2$ & $u_{1, 4} = z$ & $\mathcal R_{xz}$ invariant, $\mathcal R_{yz}$ invariant\\
	$6$ & $5$  eigenfunctions  &
	\end{tabular}\\

We are interested in the Laplace operator on $S^2_{\varphi_0}$ and not $S^2$. Roughly speaking, since $S^2_{\varphi_0}$ is just $S^2$ with very small neighborhoods removed, we expect the eigenvalues and eigenfunctions of the Laplace operator for Dirichlet problems on $S^2_{\varphi_0}$ and on $S^2$  not to differ much. The result is stated more precisely in the Appendix  B of \cite{kapouleas;constant-mean-curvature}. Combined with the symmetries of the Laplace operator, this implies that, for given $\eps>0$, there is a $\varphi_0$ small enough for which  there are eigenfunctions $u_{2, 2}, u_{2, 3}$ and $u_{2, 4}$ satisfying\\
	\begin{tabular}{c|c|c}
	Eigenvalue & Eigenfunction & Properties of eigenfunction\\ \hline
	$\lambda_{2, 1} =0$ &$ u_{2, 1}$ &  \\
	$1/2<\lambda_{2, 2} <5/2$ & $u_{2, 2} $ &$\mathcal R_{xz}$ antivariant \\
	$1/2<\lambda_{2, 3} <5/2$ & $u_{2, 3} $ & $\mathcal R_{xz}$ invariant, $\mathcal R_{yz}$ antivariant\\
	$1/2<\lambda_{2, 4} <5/2$ & $u_{2, 4} $ & $\mathcal R_{xz}$ invariant, $\mathcal R_{yz}$ invariant\\
	$13/2 < \lambda_{2, k}$ ,  $ k\geq 5 $& &
	\end{tabular}\\
	
and $\| u_{1,j} - u_{2,j} \|_{L^2(S^2_{\varphi})} \leq \eps $,  for $j=2, 3, 4$ and $0< \varphi < \varphi_0$.

\begin{theorem}
\label{thm:trivial-kernel}
There exists a $C_1$ such that for all $C_0 >C_1$, the operator $\Delta_g +|A|^2 :W^{2, p}_{Sym}(\Sigma_0^{\ast}) \cap W^{1, p}_0(\Sigma_0^{\ast}) \to L^p_{Sym}(\Sigma_0^{\ast})$ has a trivial kernel, or equivalently, there exists a $\varphi_1$ so that, for all $\varphi_0 < \varphi_1$,
 the operator $ \Delta_{g_{S^2}} +2 : W^{2, p}_{Sym}(S^2_{\varphi_0}) \cap W^{1, p}_0(S^2_{\varphi_0}) \to L^p_{Sym}(S^2_{\varphi_0}) $ has a trivial kernel. 

\end{theorem}

\begin{proof}
Suppose it is not true. Then there exists a sequence of real numbers $\varphi_j \to 0$, $j \in \mathbf{N}$ and a family of functions $u_j \in W^{2, p}_{Sym}(S^2_{\varphi_j}) \cap W^{1, p}_0(S^2_{\varphi_j}) $ so that 
	\[
	\Delta_{g_{S^2}} u_j +2 u_j = 0 \textrm{ in } S^2_{\varphi_j}, \quad \lVert u_j \rVert_{L^2(S^2_{\varphi_j})} = 4 \pi /3.
	\]
Denote by $u =u_{1,3} $ the $\mathcal R_{xz}$ invariant $\mathcal R_{yz}$ antivariant eigenfunction corresponding to the eigenvalue $2$ in $S^2$:
	\[ 
	\Delta_{g_{S^2}} u+2u = 0 \textrm{ in } S^2, \quad \lVert u \rVert_{L^2(S^2)} = 4 \pi/3.
	\]
Note that with this normalization, $u(x, y, z) =x$ for $ (x, y, z) \in S^2$. From the discussion above and \cite{kapouleas;constant-mean-curvature}, we have $u_j \to u$ as $j \to \infty$ in $L^2(K)$ for any $K$ compact, $K \subset (S^2 \setminus \{(\pm 1, 0, 0), (0, 0, \pm 1)\})$. 

\begin{claim}
\label{claim:C1-convergence}
Moreover,  $u_j \to u$ in $C^1(K)$ for any $K$ compact, $K \subset (S^2 \setminus \{(\pm 1, 0, 0), (0, 0, \pm 1)\})$. 
\end{claim}
\begin{proof}[Proof of Claim] The function $u -u_j$ satisfies 
$
	\Delta_{g_{S^2}} (u -u_j ) = 0  \textrm { in } S^2_{\varphi_j}
$. The claim is a direct consequence of interior estimates from standard elliptic theory, Sobolev inequalities and the convergence in $L^2$ on compact subsets. 
\end{proof}
Consider
	\begin{equation}
	\label{eq:sum-integrals}
		\begin{split}
	 	0 &= \int_{S^2_{\varphi_j}} u_j (\Delta_{g_{S^2}} u +2u ) - u (\Delta_{g_{S^2}}  u_j+2 u_j) \\
			&= \int_{S^2_{\varphi_j}} u_j \Delta_{g_{S^2}} u - u \Delta_{g_{S^2}} u_j =  \int_{\pd S^2_{\varphi_j}} (-u_j \nabla_{\vec{n}} u +u \nabla_{\vec{n}} u_j ) ds\\
			&= 2 \int_{c_3} u \nabla_{\vec{n}} u_j ds + \int_{c_1} u \nabla_{\vec{n}} u_j ds + \int_{c_2} u \nabla_{\vec{n}} u_j ds,
		\end{split}
	\end{equation}
where $\vec{n}$ is a the unit vector normal to $\pd S^2_{\varphi_j}$ pointing towards $S^2_{\varphi_j}$ and $c_1, c_2$ and $c_3$ are circles of (geodesic) radius $\varphi_j$ centered at $(0, 0, 1)$, $(0, 0, -1)$ and $(1, 0, 0)$ respectively. We will show that the right hand side of \eqref{eq:sum-integrals} is positive in order to obtain a contradiction. The main contribution comes from the integral over $c_3$. We bound $u$ and $\nabla u$ from below near $c_3$ with the help of a maximum principle (Lemma \ref{lem:modified-maxprinciple}) and a subsolution $w_j$. The integrals on the curves $c_1$ and $c_2$ are negligible because the eigenfunctions $u_j$ are close to $u(x, y, z) = x$, which is small near $(0, 0, \pm 1)$.  We use the convergence  in $C^1$ from Claim \ref{claim:C1-convergence} to show that the integrals on $c_1$ and $c_2$ are nonnegative. \\

We first examine $\int_{c_3} u \nabla_{\vec{n}} u_j ds$ and show that $\nabla_{\vec{n}} u_j >0$ on $c_3$. Note that the coefficient of the linear term is $2$, which has the wrong sign to apply the maximum principle directly, so we need the following lemma.

\begin{lemma}
\label{lem:modified-maxprinciple}
Let $\Omega$ be a bounded domain in $S^2$ and suppose that $u \in C^0(\bar{\Omega})$ is a supersolution of $\Delta_{g_{S^2}} +2$, i.e.
	\[
	L(u) :=\Delta_{g_{S^2}} u +2u \leq 0 \textrm{ in } \Omega
	\]
in the weak sense. 
If there exists a positive supersolution $\zeta \in C^0(\bar{\Omega})$, then $u \geq 0$ on $\pd \Omega$ implies $u \geq 0$ in $\Omega$.
\end{lemma}

\begin{proof}
Consider the function $v=\frac{u}{\zeta}$, and define the operator $\tilde{L}(v) = L(\zeta v)$. On the one hand, 
	$
	\tilde{L}(v) = L(\zeta v)  = L(u) \leq 0;
	$
but on the other hand,
	$ 
	\tilde{L}(v) = \zeta \Delta v + 2 \nabla \zeta \cdot \nabla v + (\Delta \zeta+2 \zeta) v.
	$ 
The coefficient in front of $v$ is now non positive, so we can apply the weak maximum principle to obtain
	\[
	\inf_{\Omega} v \geq \inf_{\pd \Omega} (\min(v, 0)).
	\]
Since $\zeta$ is a positive function,  $v\geq 0$ on $\pd \Omega$ therefore
	$
	\inf_{\Omega} u \geq 0.  
	$
\end{proof}

Denote by $p$ the point $(1, 0, 0)$, and choose spherical coordinates $\varphi \in [0, \pi]$ and $\theta \in [0, 2 \pi)$, 
	\[
	\begin{cases}
	x &= \cos \varphi\\
	y &= \sin \varphi \cos \theta\\
	z &= \sin \varphi \sin \theta.
	\end{cases}
	\]
The following lemma gives us the supersolution needed to apply Lemma \ref{lem:modified-maxprinciple} in the annulus $\delta_j < \varphi < \beta$, for  $\beta$ chosen below.
\begin{lemma}
\label{lem:supersol}
Let  $\beta$ be a small number so that $-\frac{\cos \beta}{\sin \beta} + 4 \beta <0$. The function $\zeta(\varphi, \theta) = 2 \beta - \varphi$ satisfies
	\[
	\Delta \zeta + 2 \zeta \leq 0 , \quad \varphi_j < \varphi< \beta
	\]
in other words, it is a supersolution of $L:= \Delta +2$ in the region $\varphi_j < \varphi< \beta$. 
\end{lemma}

\begin{proof} 
A computation using the Laplacian in spherical coordinates shows
 	\[
	\Delta \zeta +2 \zeta = -\frac{\cos \varphi}{\sin \varphi} + 2 (2 \beta -\varphi) < -\frac{\cos \beta}{\sin \beta} + 4 \beta <0, \quad \varphi_j < \varphi<\beta,
	\]
since $\cos \varphi > \cos \beta$ and $(\sin \varphi)^{-1} > (\sin \beta)^{-1}$ for $\varphi_j < \varphi<\beta$. 
\end{proof}

Define $\eta_0 = 1- \cos \beta$. The eigenfunction $u$ is given by $u= \cos \varphi$, so
	$
	| u(\varphi, \theta) - 1 | \leq \eta_0, \varphi <\beta, 0\leq \theta < 2 \pi.
	$
From the convergence of $u_j$, there is a $J$ so that 
	$
	|u_j(\beta, \theta) -1 | < 2 \eta_0,  0 \leq \theta<2\pi, j>J.
	$
For $j>J$, we have
	\begin{gather*}
	L(u_j) = \Delta u_j + 2u_j = 0, \quad \delta_j < \varphi <\beta,\\
	1+ 2 \eta_0 >u_j > 1- 2 \eta_0,   \quad \varphi = \beta, \\
	u_j = 0 , \quad   \varphi = \varphi_j.
	\end{gather*}
We introduce an auxiliary function $w_j = \varphi - \varphi_j$.  Since
	$
	L(w_j)  \geq 0, 
	$
we have $L(u_j - \eps w_j) \leq 0$ for $\varphi_j <\varphi< \beta$ and 
	$
	u_j - \eps w_j =0$ on $\varphi =\varphi_j$.
 Let $\eps = \frac{1- 2\eta_0}{\beta - \varphi_j}$, then
	\[
	u_j -\eps w_j > 1- 2\eta_0 - \frac{1- 2\eta_0}{\beta - \varphi_j} (\beta - \varphi_j) =0 \textrm{ on }\varphi=\beta.
	\]
The conditions of Lemma \ref{lem:modified-maxprinciple} are satisfied for $u = u_j -\eps w_j$ and $\Omega = \{ (\varphi, \theta) | \varphi_j < \varphi< \beta, 0\leq \theta < 2\pi \}$. Hence, 
	$
	u_j \geq \eps w_j, \quad \varphi_j <\varphi<\beta.
	$
Since $u_j = \eps w_j = 0$ on  $\varphi = \varphi_j$ (which is the circle $c_3$) we have
	\[
	\nabla_{\vec{n}} u_j \geq \eps \nabla_{\vec{n}} w_j =\eps>0 \textrm{ on } \varphi = \varphi_j
	\]
 where $\vec{n}$ is the inward unit vector normal to $c_3$. Therefore,
	\begin{equation}
	\label{eq:int-c3}
		\begin{split}
		\int_{c_3} u \nabla_{\vec{n}} u_j ds &\geq \int_{c_3} \cos \varphi_j \frac{1-2\eta_0}{\beta-\varphi_j} ds > 4 \cos\varphi_j\  \frac{1-2\eta_0}{\beta} \varphi_j \\
			& > 4 \cos \beta\left( \frac{2 \cos \beta -1}{\beta} \right)\varphi_j.
		\end{split}
	\end{equation}

Let us  now consider $ \int_{c_1} u \nabla_{\vec{n}} u_j ds$ and $ \int_{c_2} u \nabla_{\vec{n}} u_j ds$. We will show that $\nabla_{\vec{n}} u_j \geq 0$ on $c_1 \cap \{x\geq 0\}$ and $c_2 \cap \{ x\geq 0\}$.
In this case also, choose spherical coordinates $\varphi \in [0, \pi]$, $\theta \in [0, 2 \pi)$,
	\[
	\begin{cases}
	x &= \sin \varphi \cos \theta\\
	y &= \sin \varphi \sin \theta\\
	z &= \cos \varphi .
	\end{cases}
	\]
and let $\zeta$ and $\beta$ be as in Lemma \ref{lem:supersol}. By definition and symmetries, we know that for all $j$, 
	\[ 
	u_j =0 \textrm{ on } \{\varphi = \varphi_j\} \cup \{ \theta = \pm \pi/2\}.
	\]
Let $K\subset S^2\setminus \{(\pm1, 0, 0), (0, 0, \pm 1)\}$ be a compact set containing $\varphi = \beta$. Since $u_j \to u$ in $C^1(K)$, for  $\mu>0$, there is a $J$ so that 
	\[
	u - \mu \cos \theta  \leq u_j \leq u+ \mu \cos \theta \textrm{ for } \theta \in (-\pi/2, \pi/2), \varphi=\beta,  j>J.
	\]
Replacing $u$ by its value and choosing $\mu < \sin \beta$, we get
	\[
	0  \leq u_j \leq  (\sin \beta + \mu) \cos \theta, \text{ for } \theta \in (-\pi/2, \pi/2), \varphi=\beta,  j>J.
	\]
$\zeta$ is a supersolution of $\Delta +2$ in this case also, hence Lemma \ref{lem:modified-maxprinciple} applied to $u_j$ in the domain $\varphi_j <\varphi <\beta, -\pi/2 < \theta < \pi/2$ gives $u\geq 0$ in this domain. Since $u_j =0$ on $\varphi = \varphi_j$, 
	\[
	\nabla_{\vec{n}} u_j \geq 0, \quad \varphi = \varphi_j, \theta \in (-\pi/2, \pi/2).
	\]
The function $u = \sin \varphi \cos \theta $ is also non negative on $c_1 \cap \{ \theta \in (-\pi/2, \pi/2)$,  therefore
	\[
	\int_{c_1 \cap \{ -\pi/2< \theta < \pi/2\}} u \nabla_{\vec{n}} u_j ds \geq 0.
	\]
Using symmetries of $u$ and $u_j$, we have
	\begin{equation}
	\label{eq:int-c1}
	\int_{c_1} u \nabla_{\vec{n}}u_j ds = 2\int_{c_1  \cap \{ -\pi/2< \theta < \pi/2\}}u \nabla_{\vec{n}} u_j ds \geq 0.
	\end{equation}
The integral over $c_2$ is treated in an analogous way, and
	\begin{equation}
	\label{eq:int-c2}
	\int_{c_2} u \nabla_{\vec{n}}u_j ds  \geq 0.
	\end{equation}
Equations \eqref{eq:sum-integrals}, \eqref{eq:int-c3}, \eqref{eq:int-c1} and \eqref{eq:int-c2} lead to a contradiction. Therefore, there exists a $\varphi_1$ so that, for all $\varphi_0 < \varphi_1$, the operator $\Delta_{g_{S^2}} +2:W^{2, p}_{Sym}(S^2_{\varphi_0}) \cap W^{1, p}(S^2_{\varphi_0}) \to L^p(S^2_{\varphi_0})$ has a trivial kernel.
\end{proof}

\subsection{ Proof of Theorem \ref{thm:G-diffeomorphism}} 
Recall that 
	\begin{eqnarray*}
	DG(0, 0): W^{2, p}_{Sym}(\Sigma_0^{\ast}) \times \R & \to &L^p_{Sym}(\Sigma_0^{\ast}) \times \mathcal{T} \times \R\\
	(u, \sigma) &\mapsto & (\Delta u + |A|^2 u + \sigma v, [u], \sigma).
	\end{eqnarray*}
	
In what follows, all the operators and maps are considered from $W^{2, p}_{Sym}(\Sigma_0^{\ast}) \cap W^{1, p}_0(\Sigma_0^{\ast}) \to L^p_{Sym}(\Sigma_0^{\ast})$, unless mentioned otherwise. \\

The manifold $\Sigma_0^{\ast}$ is bounded and smooth, so the map $u \mapsto |A|^2 u$ is compact. The operator $\Delta_g$ is Fredholm of index $0$ since $\Delta_g = \frac{|A|^2}{2} \Delta_{\nu^{\ast} g_{S^2}}$ and since the Laplace operator on the sphere is invertible  from $W^{2, p}_{Sym}(S^2_{\varphi_0}) \cap W^{1, p}_0(S^2_{\varphi_0})$ to $L^p_{Sym}(S^2_{\varphi_0})$ for $\varphi_0 \leq 1$. Thus, 	$ \Delta_g +|A|^2$ is Fredholm of index $0$.

Let $C_1$ be given by Theorem \ref{thm:trivial-kernel}. For a choice of $C_0>C_1$, Theorem \ref{thm:trivial-kernel} and the above discussion imply that 
	$
	\Delta_g +|A|^2$ is invertible. Hence, for $v \in L^p_{Sym}(\Sigma_0^{\ast})$ and $\sigma \in \R$, 
	$
	u \mapsto \Delta_g u+|A|^2 u+\sigma v $ is invertible from $   W^{2, p}_{Sym}(\Sigma_0^{\ast}) \cap W^{1, p}_0(\Sigma_0^{\ast}) \to L^p_{Sym}(\Sigma_0^{\ast})$. We derive from this that the derivative $DG(0, 0): W^{2, p}_{Sym}(\Sigma_0^{\ast}) \times \R  \to L^p_{Sym}(\Sigma_0^{\ast}) \times \mathcal{T} \times \R$ is an isomorphism.

 We can apply the inverse function theorem for Banach spaces to obtain the existence of  a constant $\eta_0>0$, a neighborhood $U_{Sym} \subset W^{2, p}_{Sym}(\Sigma_0)$ and a neighborhood $V \subset L^p_{Sym}(\Sigma_0) \times \mathcal{T} \times (-\eta_0, \eta_0)$ such that the map $G: U_{Sym} \times (-\eta_0, \eta_0)  \to V$ is a diffeomorphism. \qed

\section{Corollaries}
\label{ssec:corollaries}
Theorem \ref{thm:existence-of-h} is an immediate corollary of Theorem \ref{thm:G-diffeomorphism}. Let us recall it:\\
\emph{There exists a constant $C_1$ so that, for any $C_0>C_1$, there are constants $\eta_0>0$ and $\delta_0>0$ such that for any $\tau\in (-\eta_0, \eta_0)$ and  function $f \in W^{2, p}_{Sym}(\Sigma_0)$ satisfying $\lVert f \rVert_{W^{2,p}(\Sigma_0)} \leq \delta_0$,
there exists a function $h \in W^{2, p}_{Sym}(\Sigma_0)$ satisfying}
	\begin{gather*}
	F(h, \tau) = H\htau+\tau e_1\cdot \nu\htau+\tau^2 X\htau \cdot \nu\htau= 0, \\
	 h=f \textrm{ on } \pd \Sigma_0 \cap (\{ x=\pm C_0\} \cup \{ z= \pm C_0\}).\end{gather*}
\emph{Moreover, $\eta_0$ and $\delta_0$ can be chosen small enough so that  $\Sigma_{h, \tau}$, the graph of $h$ over $\Sigma_{\tau}$, $\tau\in (-\eta_0, \eta_0)$, is embedded.
}\\

We wish to go back to the smaller scale, so that the surfaces considered are asymptotic to cylinders of radius $1$ instead of $1/\tau$. A few definitions are needed to state the next corollary.

For $N \in \mathbf{N}$ and a large constant $C>0$, denote by $\Omega_C^N$ the domain 
	\[
	\Omega_C^N=\{(x, y, z) \in \R^3 |\ \  |x| < C, \frac{-3N\pi}{2}<y<\frac{N\pi}{2}, |z| < C\},
	\]
and by $S_C^N$ the piece of Scherk surface
	$
	S_C^N = S \cap \Omega_C^N.
	$
For $\tau \in (-\eta, \eta)$, let $\Phi_{\tau}$ and $\mathcal H_{\tau}$ be the bending maps and the scaling maps defined by equations \eqref{eqdef:bendingmap} and \eqref{eqdef:scalingmap} respectively. The surface $\tilde \Sigma^C_{1/N} = \mathcal H_{1/N} ( \Phi_{1/N}(S_C^N))$ is a smooth surface with $N$ handles. For $x$ close to $\pm C$, $\tilde \Sigma^C_{1/N}$ is close to a plane; for $z$ close to $\pm C$,  $\tilde \Sigma^C_{1/N}$ is close to a cylinder of radius $1$. Moreover, $\tilde \Sigma^C_{1/N}$ has the following symmetries:\\
- it is  invariant with respect to the rotation of $180$ degrees about the $x$-axis and \\
- it is invariant with respect to symmetries across the planes $\theta = \pi/2N + k\pi/N$, $k =1, \ldots, 2N$, where we used cylindrical coordinates.

\begin{theorem}
\label{thm:end}
There is a constant $C_1$ so that, for any $C>C_1$, there exist an integer  $N_0>0$ and a constant $\delta_0>0$ with the following properties:\\
for every integer $N>N_0$ and for every  function $\tilde f \in W^{2, p}(\tilde \Sigma^C_{1/N})$, with  $\lVert \tilde f \rVert_{W^{2,p}(\tilde \Sigma^C_{1/N})} \leq \delta_0/(2 N^2)$ and whose graph  over $\tilde \Sigma^C_{1/N}$ satisfies the above symmetries, there exists a function $\tilde h \in W^{2,p}(\tilde \Sigma^C_{1/N})$ such that 
	\begin{gather*}
	\text{ the graph of $\tilde h$ over $\tilde \Sigma^C_{1/N}$ satisfies $H+X \cdot \nu =0$,}\\
	\text{ the  graph of $\tilde h$ over $\tilde \Sigma^C_{1/N}$ enjoys the above symmetries,}\\
	 \tilde h=\tilde f \textrm{ on } \pd \tilde \Sigma^C_{1/N}.
	\end{gather*}
Moreover, we can choose $N_0$ and $\delta_0$ so that the graph of $h$ over $\tilde \Sigma^C_{1/N}$ is an embedded surface.
\end{theorem}
\begin{proof} Lemma \ref{lem:scaling-HplusXnu} relates the equations $F(h, \tau) =0$ and $H+X\cdot \nu=0$, and the invariance with respect to the planes $\theta=\pi/2N +k\pi/N$ lets us reconstruct a complete surface from the  piece $\mathcal H_{1/N} (\Sigma _{h, 1/N})$.

Fix $C>C_1$. From now on, we will not mention the dependence in $C$ and simply write $\tilde \Sigma_{1/N}$ for $\tilde \Sigma_{1/N}^C$. 

The surface $\tilde \Sigma_{1/N}$ contains $N$  pieces of surface identical to $\mathcal H_{1/N} (\Sigma_{1/N})$. From the symmetries of $f$, we have  $\| \tilde f \|_{W^{2,p}(H_{1/N}(\Sigma_{1/N})) }\leq \delta_0/(2 N^3)$. Therefore $f = N \tilde f \circ \mathcal H_{1/N}$ considered as a function on $\Sigma_{1/N}$ satisfies
	\[
	 \| f \|_{W^{2,p}(\Sigma_{1/N})} \leq \delta_0/2.
	 \]
Moreover, $f$ enjoys the symmetries of Theorem \ref{thm:existence-of-h}. The metrics $g_{\tau}, \tau \in (-\eta, \eta)$ induce equivalent norms on the manifold $\Sigma_0$  since the bending map $\Phi$ from Lemma \ref{lem:Phitaudifferentiable} is uniformly continuous from $(-\eta, \eta) \to C^5(\bar \Omega, \R^3)$. Therefore, we can choose $\eta$ small enough so that 
	\[
	\| f \|_{W^{2, p}(\Sigma_0)} \leq 2 \| f \|_{W^{2, p} (\Sigma_{\tau})}, \quad \tau \in (-\eta, \eta).
	\]
For $N>N_0> \max(\eta, \eta_0^{-1})$ with $\eta_0$ and $\delta_0$ as in Theorem \ref{thm:existence-of-h}, there exists a function $h$ so that $F(h, 1/N)=0$ and  $\Sigma_{h,1/N}$ is embedded. Lemma \ref{lem:scaling-HplusXnu} implies that $\mathcal H_{1/N} (\Sigma_{h, 1/N})$ satisfies the self-similar equation \eqref{eq:self-shrinker}. We can define $\tilde h$ on $\mathcal H_{1/N} (\Sigma_{1/N})$ by $\tilde h = N^{-1} h \circ (\mathcal H_{1/N})^{-1}$ so that the graph of $\tilde h$ over $\mathcal H_{1/N} (\Sigma_{1/N})$ coincides with  $\mathcal H_{1/N} (\Sigma_{h, 1/N})$.  We then use invariance with respect to reflections across  the planes $ \theta=\pi/2N +k\pi/N$ to extend the function $\tilde h$ to the whole surface $\tilde \Sigma_{1/N}$.
\end{proof}
\bibliographystyle{siam}
\bibliography{thesis}

\end{document}